\documentclass[journal]{IEEEtran}

\usepackage{amsmath}
\usepackage{amssymb}
\usepackage{amsthm}
\usepackage{booktabs}
\usepackage{graphicx}
\usepackage{cite}
\usepackage{algorithm}
\usepackage{algpseudocode}
\usepackage[normalem]{ulem}
\usepackage{float}
\usepackage{verbatim}
\usepackage{threeparttable}
\begin{document}

\title{Robust Neural Policy Distillation of Long-Horizon FCS-MPC for Flying-Capacitor Three-Level Boost Converters}
\author{
Jinjian Sheng, Kazumune Hashimoto, Shuang Zhao, Mahdieh S. Sadabadi
\thanks{Jinjian Sheng and Kazumune Hashimoto are with The University of Osaka.}
\thanks{Shuang Zhao is with Hefei University of Technology.}
\thanks{Mahdieh S. Sadabadi is with The University of Manchester.}
\thanks{Corresponding author: Kazumune Hashimoto (hashimoto@eei.eng.osaka-u.ac.jp).}
}
\maketitle

\begin{abstract}
Long-horizon finite-control-set model predictive control (FCS-MPC) can improve transient regulation and flying-capacitor balancing in flying-capacitor three-level boost converters (FC-TLBCs). However, searching over switching sequences becomes computationally expensive at high switching frequencies. We train a feedforward neural network to imitate an $N$-step FCS-MPC expert computed with beam search. To improve robustness, expert trajectories are generated under randomized input voltage, load resistance, and component parameters, and a disagreement-based DAgger variant is used to relabel on-policy states where the student and expert disagree. In simulation, the learned policy maintains stable voltage regulation and capacitor balancing under nominal conditions, operating-point changes, and perturbations of several physical parameters. We demonstrate the effectiveness of our approach by reducing the computational burden. 
We also demonstrate transfer to an NPC-type three-level buck converter, where initializing from the FC-TLBC network improves sample efficiency compared with training from scratch.
\end{abstract}


\IEEEpeerreviewmaketitle

\section{Introduction}
\IEEEPARstart{F}{lying}-capacitor multilevel converters are attractive because the flying-capacitor branch enables multilevel switching operation, reduces device voltage stress, and introduces additional control freedom for capacitor-voltage balancing \cite{defay2010direct3,vyncke2012finite38}. These benefits, however, come with strongly coupled and mode-dependent dynamics among the inductor current, flying-capacitor voltage, and output voltage. As a result, achieving fast and robust closed-loop control remains difficult, especially under input-voltage sags and load variations.

Earlier work on switched and multilevel power-converter control has explored direct control strategies, weighting-factor design, and predictive-control formulations for multivariable objectives \cite{defay2010direct3,silva2009guidelines,rodriguez2012state2,vazquez2016model1}. For flying-capacitor topologies, capacitor-voltage balancing must be maintained simultaneously with output regulation and current shaping, which increases both the control complexity and the sensitivity to modeling errors and parameter mismatch \cite{vyncke2012finite38,martin2017sensitivity2}. Longer prediction horizons can improve transient behavior and steady-state performance, but the associated search complexity grows rapidly with horizon length \cite{geyer2015perfor11,keustu2024long11}.

Therefore, finite-control-set model predictive control (FCS-MPC) is a promising framework for FC-TLBCs because it evaluates admissible switching actions directly in the switching domain and can explicitly encode current tracking and flying-capacitor balancing in a common cost function \cite{rodriguez2012state2,vazquez2016model1,aguilera2012finite,vyncke2012finite38}. Its main practical limitation is computational: as the prediction horizon increases, the online search over switching sequences becomes prohibitively expensive for high switching frequencies and resource-constrained digital platforms.

To reduce this burden, neural-network approximations of MPC/FCS-MPC have been investigated for several power-electronic systems, including inverters, flying-capacitor multilevel converters, DC--DC converters, and FPGA-oriented implementations \cite{mohamed2019neural,bakeer2022artificial,wang2021model,simonetti2023neural,novak2020supervised,xiang2023light,li2025long}. These studies show that learned policies can greatly reduce inference latency. However, many are trained mainly around nominal operating conditions or evaluated under a limited set of disturbance cases. As a result, robustness to simultaneous variations in input voltage, load, and passive-component values is still not fully characterized. Moreover, pure behavior cloning is vulnerable to covariate shift: once the learned policy deviates from the expert, the closed-loop state distribution can move into regions that are rare or absent in the offline demonstrations \cite{ross2011reduction}.

This paper addresses these limitations by distilling a long-horizon FCS-MPC expert for an FC-TLBC into a compact feedforward neural policy. The expert is implemented as an $N$-step beam-search FCS-MPC controller, and its demonstrations are generated under domain randomization over operating conditions and passive-component values \cite{tobin2017domain}. To mitigate on-policy distribution shift, we further apply a disagreement-based DAgger procedure that evaluates the expert on learner-visited states and retains only disagreement states for aggregation \cite{ross2011reduction}. In this way, the proposed framework combines long-horizon expert supervision, robustness-oriented data generation, and selective on-policy relabeling within a single MPC-to-neural distillation pipeline.

The main contributions of this paper are as follows:
\begin{itemize}
\item We develop an $N$-step beam-search FCS-MPC expert for FC-TLBC inner-loop control and distill it into a four-class feedforward neural switching policy.
\item We propose a robust data-generation and imitation-learning pipeline that combines domain randomization over operating points and passive components with selective on-policy relabeling via disagreement-based DAgger.
\item We present scenario-based simulation results showing stable regulation, current tracking, and flying-capacitor balancing under nominal conditions, operating-point variations, and perturbations in $L$, $C_f$, and $C$, while substantially reducing the online decision time relative to the expert on the same evaluation CPU.
\item We demonstrate transfer to an NPC-type three-level buck converter, where initialization from the FC-TLBC policy improves sample efficiency relative to training from scratch.
\end{itemize}
\textit{Related work.}\\
\smallskip
\noindent\textbf{Predictive Control for Switched and Multilevel Converters.}
For switched and multilevel converters, predictive control is attractive because it operates directly in the switching domain and can handle current tracking, voltage regulation, and capacitor balancing within a unified optimization framework \cite{kouro2008model3,silva2009guidelines,rodriguez2012state2,vazquez2016model1}. In the broader control-systems literature, implementation-oriented predictive and hybrid-control studies have also been reported for step-down, buck/boost, full-bridge, and boost DC--DC converters \cite{geyer2008hybrid,mariethoz2010comparison,xie2012fullbridge,kim2014stabilizing}. For flying-capacitor and related multilevel topologies, prior studies have shown that predictive formulations are particularly useful when internal capacitor-voltage balancing must be coordinated with external regulation objectives \cite{defay2010direct3,vyncke2012finite38,scoltock2015model3}. Compared with shorter-horizon or simplified predictive strategies, longer-horizon formulations can improve transient behavior and steady-state quality, but the online combinatorial search grows rapidly with the horizon length and the number of admissible switching actions \cite{aguilera2012finite,geyer2015perfor11,keustu2024long11}.

\smallskip

\noindent\textbf{Learning-Based Approximations of MPC/FCS-MPC.}
To reduce the online computational burden, neural-network approximations of MPC/FCS-MPC have been investigated for inverter systems with output filters, flying-capacitor multilevel converters, rectifiers, and DC--DC converters \cite{mohamed2019neural,bakeer2022artificial,liu2024finite,wang2021model,xiang2023light,novak2020supervised}. Compared with solving the predictive optimization problem at every sampling instant, these learned surrogates offer much lower inference latency and are therefore attractive for fast digital implementation. This line of work is also consistent with the long-standing emphasis on computational tractability and sampled-data implementation in predictive control of converter systems \cite{geyer2008hybrid,almer2013sampled,keustu2024long11}. Hardware-oriented studies have also been reported for converter families such as CHB topologies and for long-horizon data-driven control pipelines \cite{simonetti2023neural,li2025long}. However, many existing studies focus mainly on nominal-condition training or evaluate robustness only under a limited set of disturbances.

\smallskip

\noindent\textbf{Imitation Learning Under Distribution Shift and Robustness.}
Pure behavior cloning from offline expert trajectories is simple and effective, but compared with on-policy aggregation methods it is more vulnerable to covariate shift: once the learned controller deviates from the expert, the closed-loop trajectory may move into state regions that are weakly represented in the training data \cite{ross2011reduction}. Domain randomization addresses a complementary issue by broadening the training distribution over operating conditions and parameter values, thereby improving generalization to unseen scenarios \cite{tobin2017domain}. In related predictive-control work, practical issues such as sampled-data behavior, parameter variation, and performance adaptation have also been emphasized in converter applications \cite{martin2017sensitivity2,Novak2018statistical5,almer2013sampled,yang2018adaptive}. Nevertheless, their integration with long-horizon FCS-MPC distillation remains limited.

\smallskip

\noindent\textbf{Positioning of This Work.}
Compared with prior studies that typically emphasize either fast neural approximation or limited robustness evaluation, and compared with implementation-oriented predictive-control studies that do not consider neural distillation, the present work combines four elements in a single framework: a long-horizon beam-search FCS-MPC expert, domain-randomized expert data over both operating conditions and passive-component values, selective on-policy relabeling via disagreement-based DAgger, and scenario-based validation on an FC-TLBC under input-voltage, load, and parameter perturbations. This combination is intended to preserve the benefits of long-horizon predictive control while reducing online computational cost and improving robustness to closed-loop distribution shift.


\begin{figure}
    \centering
    \includegraphics[width=1\linewidth]{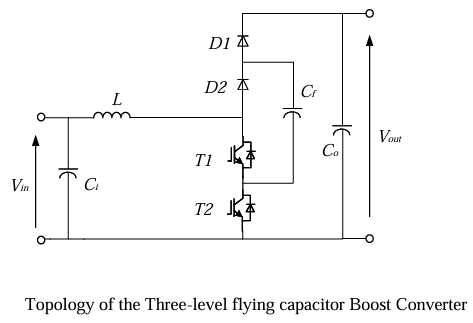}
    \caption{FC-TLBC Topology}
    \label{fig:topology}
\end{figure}

\section{Problem Formulation and Converter Model}

\subsection{Problem Setup and Feasible Switching Modes}

We consider inner-loop control of the flying-capacitor three-level boost converter (FC-TLBC) shown in Fig.~1. The overall closed-loop objective is to regulate the output voltage $v_o$ to a prescribed reference $v_o^{\star}$ while maintaining the flying-capacitor voltage as:
\begin{equation}
V_{Cf}^{\star}=\frac{v_o^{\star}}{2}.
\label{eq:VCf_ref}
\end{equation}
Following a cascaded design, an outer voltage controller generates the inductor-current reference $i_{\mathrm{ref}}$, and the inner-loop controller selects one admissible switching mode at each sampling instant.

The converter state and measurable exogenous input are defined as
\begin{equation}
\begin{aligned}
\mathbf{x}(t)=
\begin{bmatrix}
 i_L(t)\\
 v_{Cf}(t)\\
 v_o(t)
\end{bmatrix},\ 
\mathbf{w}(t)=
\begin{bmatrix}
 V_{\mathrm{in}}(t)\\
 i_o(t)
\end{bmatrix}.
\end{aligned}
\label{eq:state_input}
\end{equation}
where $i_L$ is the inductor current, $v_{Cf}$ is the flying-capacitor voltage, $v_o$ is the output voltage, $V_{\mathrm{in}}$ is the input voltage, and $i_o$ is the output current. Notice here that $\mathbf{w}(t)$ is not a control input; it is a measurable exogenous input used by the prediction model. The inner-loop manipulated variable is the admissible switching mode selected at each sampling instant, which will be described below.

To describe the switching behavior, we use a symbolic mode encoding $m=(S_A,S_B)$, where $S_A\in\{\mathrm{P},\mathrm{O},\mathrm{N}\}$ denotes the inductor terminal-voltage level and $S_B\in\{\mathrm{P},\mathrm{O},\mathrm{N}\}$ denotes the charging direction of the flying capacitor. This is a functional encoding of the converter mode rather than a direct listing of binary gate signals. In particular, $S_A=\mathrm{P}$, $\mathrm{O}$, and $\mathrm{N}$ correspond to positive, intermediate, and negative inductor terminal-voltage levels, respectively, while $S_B=\mathrm{P}$, $\mathrm{O}$, and $\mathrm{N}$ denote forward charging, no net charge transfer, and reverse charging of the flying capacitor.

Although the symbolic grid $\{\mathrm{P},\mathrm{O},\mathrm{N}\}^2$ contains nine combinations, topological constraints and Kirchhoff's voltage law reduce the admissible set to four feasible modes:
\begin{equation}
\mathcal{U}=\{\mathrm{OP},\mathrm{PO},\mathrm{NO},\mathrm{ON}\}.
\label{eq:admissible_modes}
\end{equation}
These feasible combinations are summarized in Table~\ref{tab:switching}.

\begin{table}[t]
\centering
\caption{Viable Switching Combinations for FC-TLBC}
\label{tab:switching}
\begin{tabular}{c|ccc}
\toprule
$S_B \backslash S_A$ & $\mathbf{P}$ & $\mathbf{O}$ & $\mathbf{N}$ \\
\midrule
$\mathbf{P}$ & \sout{PP} & OP & \sout{NP} \\
$\mathbf{O}$ & PO & \sout{OO} & NO \\
$\mathbf{N}$ & \sout{PN} & ON & \sout{NN} \\
\bottomrule
\end{tabular}
\end{table}

\subsection{Mode-Dependent Prediction Model}

For each feasible mode $m\in\mathcal{U}$, the FC-TLBC is represented by a mode-dependent affine state-space model. Using the state vector in \eqref{eq:state_input} and the measurable exogenous input $\mathbf{w}=[V_{\mathrm{in}},\ i_o]^\top$, we write
\begin{equation}
\dot{\mathbf{x}}(t)=\mathbf{A}_m\mathbf{x}(t)+\mathbf{B}\mathbf{w}(t).
\label{eq:mode_ct}
\end{equation}
The measured output current $i_o$ is used directly as an exogenous input, so the predictive model does not require an explicit load parameter. Equivalently, one may estimate $R_k\approx v_{o,k}/i_{o,k}$ when needed, but the rollout below only requires $i_o$.

We parameterize the four feasible modes using coefficients $\{a_{vo}^{(m)},a_{Cf}^{(m)},\alpha^{(m)},\beta^{(m)}\}$ such that
\begin{align}
\dot{i}_L &= \frac{1}{L}\Bigl(V_{\mathrm{in}}-a_{vo}^{(m)}v_o-a_{Cf}^{(m)}v_{Cf}\Bigr), \label{eq:mode_il}\\
\dot{v}_{Cf} &= \frac{\beta^{(m)}}{C_f} i_L, \label{eq:mode_vcf}\\
\dot{v}_o &= \frac{1}{C}\Bigl(\alpha^{(m)} i_L - i_o\Bigr). \label{eq:mode_vo}
\end{align}
Therefore, we have
\begin{equation}
\begin{aligned}
\mathbf{A}_m&=
\begin{bmatrix}
0 & -\dfrac{a_{Cf}^{(m)}}{L} & -\dfrac{a_{vo}^{(m)}}{L}\\
\dfrac{\beta^{(m)}}{C_f} & 0 & 0\\
\dfrac{\alpha^{(m)}}{C} & 0 & 0
\end{bmatrix},\\
\mathbf{B}&=
\begin{bmatrix}
\dfrac{1}{L} & 0\\
0 & 0\\
0 & -\dfrac{1}{C}
\end{bmatrix}.
\end{aligned}
\label{eq:mode_AB}
\end{equation}
The corresponding mode coefficients are listed in Table~\ref{tab:mode_coeff}.

\begin{table}[t]
\centering
\caption{Mode Coefficients for OP/PO/NO/ON Used in \eqref{eq:mode_il}--\eqref{eq:mode_vo}}
\label{tab:mode_coeff}
\begin{tabular}{c c c c c}
\toprule
Mode $m$ & $a_{vo}^{(m)}$ & $a_{Cf}^{(m)}$ & $\alpha^{(m)}$ & $\beta^{(m)}$ \\
\midrule
NO & 0 & 0  & 0 & 0 \\
PO & 1 & 0  & 1 & 0 \\
OP & 0 & 1  & 1 & 1 \\
ON & 1 & -1 & 1 & -1 \\
\bottomrule
\end{tabular}
\end{table}

Using forward Euler discretization with sampling period $T_s$, \eqref{eq:mode_ct} yields the discrete-time prediction model
\begin{equation}
\begin{aligned}
\mathbf{x}_{k+1}&=\mathbf{A}_{d,m}\mathbf{x}_k+\mathbf{B}_d\mathbf{w}_k,\\
\mathbf{A}_{d,m}&=\mathbf{I}+T_s\mathbf{A}_m,\  \mathbf{B}_d=T_s\mathbf{B}, 
\end{aligned}
\label{eq:mode_dt}
\end{equation}
which is used by the FCS-MPC expert during finite-horizon rollout. 

\subsection{Control Objective and Constraints}

The control objective is defined in a cascaded manner. At the closed-loop level, the converter should regulate the output voltage \(v_o\) to the reference \(v_o^\star\) while maintaining the flying-capacitor voltage around
\[
V_{Cf}^{\star} = \frac{v_o^\star}{2}.
\]
To achieve this, the outer voltage controller converts the output-voltage regulation task into an inductor-current reference \(i_{\mathrm{ref}}\). The inner-loop switching controller then selects one admissible mode \(m_k \in \mathcal{U}\) at each sampling instant so as to (i) make \(i_L\) track \(i_{\mathrm{ref}}\), (ii) keep \(v_{Cf}\) close to \(V_{Cf}^{\star}\), and (iii) satisfy the hard current limit.
Thus, \(v_o\) is regulated indirectly through the outer loop, whereas the inner loop acts directly on the switching mode. 

\section{Proposed MPC-to-Neural Distillation Framework}

\subsection{Overview of the Proposed Framework}

The proposed workflow consists of four stages: (i) construct a long-horizon FCS-MPC expert based on the prediction model in Section~II; (ii) generate expert demonstrations under randomized operating conditions and parameter values; (iii) train a compact feedforward neural policy to imitate the expert's switching decision; and (iv) refine the policy with selective on-policy relabeling using disagreement-based DAgger. The goal is to retain the closed-loop behavior of long-horizon predictive control while reducing the online decision cost to that of a simple classifier.

\subsection{N-Step Beam-Search FCS-MPC Expert}

At each sampling instant, the expert receives the measured information vector
\begin{equation}
\mathbf{z}_k=
\begin{bmatrix}
 i_{L,k}\\
 v_{Cf,k}\\
 v_{o,k}\\
 i_{\mathrm{ref},k}\\
 V_{\mathrm{in},k}\\
 i_{o,k}
\end{bmatrix},
\label{eq:feature_vector}
\end{equation}
which contains the plant state, the outer-loop current reference, and the measurable exogenous quantities required by the prediction model. The expert then evaluates a candidate mode sequence
\begin{equation}
\begin{aligned}
m_{k:k+N-1} &= \{m_k,m_{k+1},\ldots,m_{k+N-1}\},\\
m_{k+j} &\in \mathcal{U},\; j=0,\ldots,N-1.
\end{aligned}
\label{eq:candidate_sequence}
\end{equation}
by rolling out the mode-dependent model in \eqref{eq:mode_dt}. The associated finite-horizon cost is
\begin{multline}
J(m_{k:k+N-1})=
\sum_{n=1}^{N}
\Bigl[
\lambda_I\bigl(i_{L,k+n}-i_{\mathrm{ref},k+n}\bigr)^2 \\
+ \lambda_{Cf}\bigl(v_{Cf,k+n}-V_{Cf}^{\star}\bigr)^2
\Bigr], 
\label{eq:mpc_cost}
\end{multline}
where $\lambda_I$, $\lambda_{Cf}$ are the weight parameters. 
The optimal sequence is defined as
\begin{equation}
m_{k:k+N-1}^{\star}
=
\arg\min_{m_{k:k+N-1}} J(m_{k:k+N-1}),
\label{eq:mpc_seq}
\end{equation}
and the expert applies only the first element in receding-horizon fashion:
\begin{equation}
\pi_{\mathrm{MPC}}(\mathbf{z}_k)=m_k^{\star}.
\label{eq:mpc_policy}
\end{equation}

A naive exhaustive search would enumerate all length-\(N\) mode sequences in \(\mathcal{U}^N\), which requires \(|\mathcal{U}|^N\) complete-sequence evaluations at each sampling instant. For the FC-TLBC considered here, \(|\mathcal{U}|=4\), so exhaustive search already involves \(4^N\) candidate sequences (e.g., \(1024\) when \(N=5\)). To reduce this burden, we employ beam search, which grows the search tree stage by stage rather than enumerating all complete sequences. At depth \(\ell\), each retained partial sequence \(m_{k:k+\ell-1}\) is expanded by all admissible next modes in \(\mathcal{U}\), the cumulative cost of the resulting children is updated, and only the \(K\) partial sequences with the lowest cumulative cost are kept for the next expansion. After the tree reaches depth \(N\), the complete sequence with the smallest cost is selected, and only its first mode is applied in receding-horizon fashion. The number of candidate expansions is therefore on the order of \(K|\mathcal{U}|N\), which is much smaller than \(|\mathcal{U}|^N\) when \(K \ll |\mathcal{U}|^{N-1}\). The price paid for this reduction is approximate optimality, since a branch discarded at an intermediate depth cannot be recovered later. Nevertheless, beam search preserves multi-step look-ahead while keeping the online computation manageable.

\subsection{Domain-Randomized Expert Dataset Construction}

To improve robustness to operating-point shifts and parameter mismatch, expert demonstrations are collected under randomized environments rather than under a single nominal condition. For each sampled environment, the expert policy $\pi_{\mathrm{MPC}}$ in \eqref{eq:mpc_policy} is executed in closed loop, and the resulting state--mode pairs are recorded.

We consider two sources of variability. The first is operating-condition variability, represented by changes in input voltage and load. The second is parameter variability, represented by perturbations in the passive components $L$, $C_f$, and $C$. The perturbed components are modeled as
\begin{equation}
\begin{aligned}
L'&=(1+\delta_L)L,\\
C_f'&=(1+\delta_{C_f})C_f,\\
C'&=(1+\delta_C)C,
\end{aligned}
\label{eq:dr_scale}
\end{equation}
where $\delta_L$, $\delta_{C_f}$, and $\delta_C$ are sampled from prescribed bounded distributions. Likewise, the operating conditions are generated by sampling the input voltage and load from predefined distributions. The exact numerical ranges used in the experiments are specified in Section~IV-A.

Each dataset sample consists of the measured feature vector in \eqref{eq:feature_vector} and its expert label:
\begin{equation}
\bigl(\mathbf{z}_k,\,\pi_{\mathrm{MPC}}(\mathbf{z}_k)\bigr),
\label{eq:sample_def}
\end{equation}
where $\pi_{\mathrm{MPC}}(\mathbf{z}_k)\in\mathcal{U}$. Although load conditions are randomized during data generation, the student policy does not require direct access to the load parameter. Instead, it uses the measurable output current $i_o$, which makes the learned controller deployable without load-parameter identification.

The offline dataset is assembled from three subsets: a nominal subset for basic steady-state and transient behavior, an operating-point-randomized subset for broader coverage of input-voltage and load variations, and a parameter-randomized subset for robustness to passive-component mismatch. Denoting these subsets by $\mathcal{D}_{\mathrm{nom}}$, $\mathcal{D}_{\mathrm{op}}$, and $\mathcal{D}_{\mathrm{par}}$, respectively, the combined dataset is
\begin{equation}
\mathcal{D}_{\mathrm{DR}}
=
\mathcal{D}_{\mathrm{nom}}
\cup
\mathcal{D}_{\mathrm{op}}
\cup
\mathcal{D}_{\mathrm{par}}.
\label{eq:dr_union}
\end{equation}

\subsection{Neural Policy and Supervised Distillation}

The student policy $\pi_{\mathrm{ANN}}$ is a compact feedforward classifier that maps the six-dimensional feature vector $\mathbf{z}_k$ in \eqref{eq:feature_vector} to one of the four admissible switching modes in $\mathcal{U}$. The policy definition used in the distillation process is summarized in Table~\ref{tab:impl_nn}.

\begin{table}[t]
\centering
\caption{Neural Policy Definition Used in Distillation}
\label{tab:impl_nn}
\begin{tabular}{l l}
\toprule
\textbf{Item} & \textbf{Setting} \\
\midrule
Network type & Feedforward fully connected classifier \\
Input features & $(i_L,\ v_{Cf},\ v_o,\ i_{\mathrm{ref}},\ V_{\mathrm{in}},\ i_o)$ \\
Output & 4 admissible switching modes \\
Output layer & Softmax classifier \\
Loss function & Class-weighted cross-entropy \\
Domain randomization & $V_{\mathrm{in}},\ R,\ L,\ C_f,\ C$ \\
On-policy correction & Disagreement-based DAgger relabeling \\
\bottomrule
\end{tabular}
\end{table}

Let $\hat{\mathbf{y}}\in\mathbb{R}^4$ denote the output class-probability vector. A representative simple feedforward policy can be written as
\begin{equation}
\hat{\mathbf{y}}=\mathrm{softmax}\!\left(\mathbf{W}_3\,\sigma\!\left(\mathbf{W}_2\,\sigma\!\left(\mathbf{W}_1\mathbf{z}+\mathbf{b}_1\right)+\mathbf{b}_2\right)+\mathbf{b}_3\right),
\label{eq:nn_softmax}
\end{equation}
where $\sigma(\cdot)$ is the activation function. The corresponding switching decision is denoted by $\pi_{\mathrm{ANN}}(\mathbf{z}_k)\in\mathcal{U}$. The student is trained by minimizing the class-weighted cross-entropy loss
\begin{equation}
\mathcal{L}(\theta)=-\sum_{c=1}^{4}\alpha_c\,y_c\log \hat{y}_c,
\label{eq:ce_loss}
\end{equation}
where $\mathbf{y}$ is the one-hot expert label and $\alpha_c$ is inversely proportional to the class frequency of class $c$.

Training on $\mathcal{D}_{\mathrm{DR}}$ alone corresponds to standard behavior cloning. While domain randomization broadens coverage over operating conditions and parameter values, it does not guarantee coverage of the state distribution actually visited by the learned policy during closed-loop execution. This motivates the on-policy refinement step described next.

\subsection{Disagreement-Based DAgger Refinement}
Behavior cloning on the offline dataset $\mathcal{D}_{\mathrm{DR}}$ provides an initial student policy, but it remains vulnerable to covariate shift. Once the learned controller deviates from the expert, the closed-loop trajectory may move into state regions that are weakly represented in the offline demonstrations, and the resulting errors can accumulate over time. To mitigate this effect, we adopt DAgger \cite{ross2011reduction}. In standard DAgger, the current learner is rolled out in closed loop, the expert is evaluated on the learner-visited states, and those on-policy states are aggregated into the training set for iterative retraining. In this paper, we use a disagreement-filtered variant of DAgger. The iterative structure is the same as in standard DAgger, but instead of aggregating all learner-visited states, we retain only those states on which the student and expert choose different switching modes. This focuses refinement on weakly cloned or failure-prone regions of the state space while keeping the additional dataset compact.

Let $\pi_{\mathrm{ANN}}^{(i)}$ denote the student policy at DAgger iteration $i$, and initialize the aggregated dataset as
\begin{equation}
\mathcal{D}_{\mathrm{aug}}^{(0)} = \mathcal{D}_{\mathrm{DR}}.
\label{eq:dagger_init}
\end{equation}
During rollout of $\pi_{\mathrm{ANN}}^{(i)}$, we evaluate the expert on the learner-visited states and define the mismatch set as
\begin{equation}
\begin{aligned}
\mathcal{D}_{\mathrm{mist}}^{(i)}
=
\Bigl\{
\bigl(\mathbf{z}_k,\,\pi_{\mathrm{MPC}}(\mathbf{z}_k)\bigr)
\;\Big|\;
\pi_{\mathrm{ANN}}^{(i)}(\mathbf{z}_k)
\neq
\pi_{\mathrm{MPC}}(\mathbf{z}_k)
\Bigr\}.
\end{aligned}
\label{eq:dagger_mismatch}
\end{equation}
The aggregated dataset is then updated by
\begin{equation}
\mathcal{D}_{\mathrm{aug}}^{(i+1)}
=
\mathcal{D}_{\mathrm{aug}}^{(i)}
\cup
\mathcal{D}_{\mathrm{mist}}^{(i)},
\label{eq:dagger_union}
\end{equation}
and the student is fine-tuned on $\mathcal{D}_{\mathrm{aug}}^{(i+1)}$. Repeating this procedure reduces on-policy distribution mismatch and improves robustness when the learner induces state trajectories that differ from those in the original offline dataset.

Algorithm~\ref{alg:dagger} summarizes one refinement cycle. Starting from the offline-trained student, each iteration performs closed-loop rollouts with the current student policy, evaluates the expert at the visited states, stores only disagreement states, augments the aggregated dataset, and fine-tunes the student on the updated dataset. Relative to standard DAgger, the key difference is therefore the filtering rule before aggregation: only disagreement samples are retained.

\begin{algorithm}[t]
\caption{Disagreement-Based DAgger Refinement}
\label{alg:dagger}
\begin{algorithmic}[1]
\State \textbf{Input:} offline dataset $\mathcal{D}_{\mathrm{DR}}$, expert $\pi_{\mathrm{MPC}}$
\State \textbf{Input:} initial student $\pi_{\mathrm{ANN}}^{(0)}$, number of iterations $I$
\State $\mathcal{D}_{\mathrm{aug}}^{(0)} \gets \mathcal{D}_{\mathrm{DR}}$
\For{$i = 0,1,\dots,I-1$}
    \State $\mathcal{D}_{\mathrm{mist}}^{(i)} \gets \emptyset$
    \For{each closed-loop rollout episode}
        \State Initialize the converter state
        \For{each time step in the rollout horizon}
            \State Observe $\mathbf{z}_k$ from the learner-induced trajectory
            \State Compute $\pi_{\mathrm{ANN}}^{(i)}(\mathbf{z}_k)$ and $\pi_{\mathrm{MPC}}(\mathbf{z}_k)$
            \If{$\pi_{\mathrm{ANN}}^{(i)}(\mathbf{z}_k) \neq \pi_{\mathrm{MPC}}(\mathbf{z}_k)$}
                \State Add $\bigl(\mathbf{z}_k, \pi_{\mathrm{MPC}}(\mathbf{z}_k)\bigr)$ to $\mathcal{D}_{\mathrm{mist}}^{(i)}$
            \EndIf
            \State Apply $\pi_{\mathrm{ANN}}^{(i)}(\mathbf{z}_k)$ to the closed loop
        \EndFor
    \EndFor
    \State $\mathcal{D}_{\mathrm{aug}}^{(i+1)} \gets \mathcal{D}_{\mathrm{aug}}^{(i)} \cup \mathcal{D}_{\mathrm{mist}}^{(i)}$
    \State Fine-tune the student on $\mathcal{D}_{\mathrm{aug}}^{(i+1)}$ to obtain $\pi_{\mathrm{ANN}}^{(i+1)}$
\EndFor
\State \textbf{Output:} refined student policy $\pi_{\mathrm{ANN}}^{(I)}$
\end{algorithmic}
\end{algorithm}

\begin{table*}[t]
\centering
\caption{Overview of Experimental Modules and Their Roles}
\label{tab:exp_overview}
\begin{tabular}{l l}
\toprule
Module & Main Purpose \\
\midrule
Basic Experiments & Compare ANN vs. MPC in S1--S3 \\
Ablation Study & Quantify roles of DR, Disagreement-Based DAgger, and expert supervision \\
Sensitivity & Robustness to DR range and Disagreement-Based DAgger budget \\
Transfer Learning & Cross-topology generalization (Buck-3L) \\
\bottomrule
\end{tabular}
\end{table*}

\section{Simulation Setup and Validation}

\subsection{Experimental Setup and Common Settings}

To evaluate the proposed framework, expert-data generation, policy distillation, and network training are conducted offline in Python/PyTorch on an Apple M3 Max CPU.
The trained ANN controller is then deployed in a Simulink model of the FC-TLBC for closed-loop validation, and the same CPU is used for the runtime comparison reported in Section~IV-B1.
The nominal converter parameters, expert-controller configuration, and ANN training settings used in the main FC-TLBC experiments are summarized in Tables~\ref{tab:converter_params}, \ref{tab:mpc_params}, and \ref{tab:ann_cfg}, respectively.

\begin{table}[t]
\centering
\caption{Nominal Parameters of the FC-TLBC and Control-Update Settings}
\label{tab:converter_params}
\begin{tabular}{l c}
\toprule
\textbf{Item} & \textbf{Value} \\
\midrule
Input voltage (nominal) \(V_{\mathrm{in}}\) & \(120~\mathrm{V}\) \\
Output reference \(v_o^\star\) & \(180~\mathrm{V}\) \\
Inductor \(L\) & \(1~\mathrm{mH}\) \\
Flying capacitor \(C_f\) & \(50~\mu\mathrm{F}\) \\
Output capacitor \(C\) & \(125~\mu\mathrm{F}\) \\
Load resistance (nominal) \(R\) & \(36~\Omega\) \\
Flying-capacitor reference \(V_{Cf}^\star\) & \(v_o^\star/2 = 90~\mathrm{V}\) \\
Control-update period \(T_{s}\) & \(20~\mu\mathrm{s}\) \\
\bottomrule
\end{tabular}
\end{table}

\begin{table}[t]
\centering
\caption{Configuration of the \(N\)-Step Beam-Search FCS-MPC Expert}
\label{tab:mpc_params}
\small
\setlength{\tabcolsep}{4pt}
\renewcommand{\arraystretch}{1.05}
\begin{tabular}{@{}p{0.56\columnwidth}p{0.2\columnwidth}@{}}
\toprule
\textbf{Item} & \textbf{Value} \\
\midrule
Action set size \(|\mathcal{U}|\) & 4\\
Prediction horizon \(N\) & 5 \\
Beam width \(K\) & 15 \\
Current tracking weight \(\lambda_I\) & 1.0 \\
Flying-capacitor voltage weight \(\lambda_{Cf}\) & 0.007 \\
\bottomrule
\end{tabular}
\end{table}

\begin{table}[t]
\centering
\caption{ANN Architecture and Training Hyperparameters Used in the Basic Experiments}
\label{tab:ann_cfg}
\small
\begin{tabular}{lc}
\toprule
Item & Value \\
\midrule
Input dimension & 6 \\
Hidden layers & 1 \\
Hidden units per hidden layer & 128 \\
Output classes & 4 \\
Activation function & ReLU \\
Output layer & Softmax \\
Optimizer & Adam \\
Learning rate & \(1\times10^{-4}\) \\
Batch size & 2048 \\
Weight decay & None \\
Offline-training epochs & 260 \\
DAgger fine-tuning epochs & 280 \\
Numerical precision & float32 \\
\bottomrule
\end{tabular}
\end{table}
The evaluation is organized into the four modules summarized in Table~\ref{tab:exp_overview}.
Basic experiments compare the distilled ANN policy against the \(N\)-step FCS-MPC expert under Scenarios~S1--S3.
Ablation removes DR or Disagreement-Based DAgger to isolate their effects under fixed test trajectories.
Sensitivity sweeps the DR range and the Disagreement-Based DAgger mismatch-sample budget to assess training robustness.
Finally, transfer learning evaluates whether features learned on FC-TLBC accelerate training and improve closed-loop behavior on an NPC-type three-level buck converter (Buck-3L).

\paragraph{ANN architecture and training details.}
The student policy is implemented as a fully connected feedforward classifier with six input features (see \eqref{eq:feature_vector}),
one hidden layer of 128 units, and a four-class softmax output corresponding to the four admissible switching modes.
The hidden layer uses the ReLU activation function. The base offline training is run for 260 epochs, followed by 280 epochs of fine-tuning after disagreement-based DAgger aggregation.
Unless otherwise stated, all reported FC-TLBC results use this same architecture and hyperparameter setting.

In the numerical experiments, the three dataset subsets introduced in Section~III-C are instantiated as Scenarios~S1--S3.
Scenario~S1 uses representative nominal step responses to capture baseline transient and steady-state behavior.
Scenario~S2 broadens the operating conditions by sampling
\begin{equation}
\begin{aligned}
V_{\mathrm{in}} &\sim \mathcal{U}(80,140)~\mathrm{V},\\
R &\sim \mathcal{U}(10,100)~\Omega,
\end{aligned}
\label{eq:scenario_s2}
\end{equation}
while Scenario~S3 further introduces passive-component perturbations according to \eqref{eq:dr_scale} with independent perturbations
\begin{equation}
\delta_L,\ \delta_{C_f},\ \delta_C \sim \mathcal{U}(-\rho,\rho),
\label{eq:scenario_s3_rho}
\end{equation}
where \(\rho\in(0,1)\) denotes the relative randomization intensity.
In this experiments, we use \(\rho=0.3\).

Algorithm~\ref{alg:basic_fc_tlbc} summarizes the data-generation, offline training, on-policy refinement, and evaluation workflow for Scenarios~S1--S3.

\begin{algorithm}[t]
\caption{Experimental Pipeline for FC-TLBC (Scenarios S1--S3)}
\label{alg:basic_fc_tlbc}
\begin{algorithmic}[1]
\State Set the random seed and simulation parameters
\State Initialize the \(N\)-step beam-search FCS-MPC expert
\For{scenario \(s \in \{S1, S2, S3\}\)}
    \State Simulate the FC-TLBC under scenario-specific \(V_{\mathrm{in}}(t)\), \(R(t)\), and parameter settings
    \State Collect expert-labeled pairs \(\bigl(\mathbf{z}_k, u_{\mathrm{MPC},k}\bigr)\) to form \(\mathcal{D}_s\)
\EndFor
\State Construct \(\mathcal{D}_{\mathrm{DR}} = \mathcal{D}_{S1} \cup \mathcal{D}_{S2} \cup \mathcal{D}_{S3}\)
\State Train the ANN on \(\mathcal{D}_{\mathrm{DR}}\) using weighted cross-entropy
\State Refine the ANN with disagreement-based DAgger to obtain the final policy
\For{scenario \(s \in \{S1, S2, S3\}\)}
    \State Run closed-loop simulations with the MPC expert and the ANN policy
    \State Compute tracking, transient, energy, penalty, and switching metrics
    \State Compare the controllers under the same scenario
\EndFor
\end{algorithmic}
\end{algorithm}

\subsection{Comparative Experiments (Scenarios 1--3)}
\begin{table*}[t]
\footnotesize
\centering
\begin{threeparttable}
\caption{Consolidated closed-loop results for the FC-TLBC under Scenarios~1--3. Scenario~3 reports the ANN only, since this case is used primarily for robustness validation under plant mismatch. Detailed metric definitions are given in Appendix~\ref{sec:metric_defs}.}
\label{tab:key_results_all}
\begin{tabular}{lccccc}
\toprule
 & \multicolumn{2}{c}{Scenario 1} & \multicolumn{2}{c}{Scenario 2} & Scenario 3 \\
\cmidrule(lr){2-3} \cmidrule(lr){4-5} \cmidrule(lr){6-6}
Metric & MPC & ANN & MPC & ANN & ANN \\
\midrule
Decision time (\(\mu\mathrm{s}\)) & 342.26 & 18.30 & 342.26 & 18.30 & 18.30 \\
Runtime (s) & 17.46 & 1.42 & 16.80 & 1.45 & 1.42 \\
\midrule
\(\mathrm{MSE}_{v_{Cf}}\) & 1.803 & 1.664 & 0.952 & 0.962 & 5.170 \\
\(\mathrm{MSE}_{v_o}\) & 14.13 & 6.22 & 10.81 & 4.03 & 33.94 \\
\(\mathrm{MSE}_{i_L}\) & 0.206 & 0.096 & 0.175 & 0.076 & 0.259 \\
\midrule
\(\mathrm{Overshoot}_{v_{Cf}}\) (V) & 8.16 & 4.65 & 3.70 & 3.90 & 26.95 \\
\(\mathrm{Overshoot}_{v_o}\) (V) & 8.93 & 33.39 & 0.69 & 25.49 & 48.28 \\
\midrule
\(N_{i_L,\mathrm{viol}}\) & 0 & 0 & 0 & 0 & 0 \\
\(\mathrm{Penalty}_{\mathrm{over}}\) & 0 & 0.0002 & 0 & 0.0001 & 0.0004 \\
\(\mathrm{Penalty}_{\mathrm{sag}}\) & 0.0007 & 0.0001 & 0.0006 & 0.0001 & 0.0011 \\
\bottomrule
\end{tabular}
\begin{tablenotes}
\footnotesize
\item Energy- and switching-related quantities are omitted from the main-text table for brevity and may be retained in the Appendix if desired.
\end{tablenotes}
\end{threeparttable}
\end{table*}

Table~\ref{tab:key_results_all} shows a consistent pattern across the three scenarios. Detailed metric definitions are given in Appendix~\ref{sec:metric_defs}. 

\subsubsection{Scenario 1: Nominal Operating Condition}
\label{sec:s1}
\begin{figure}
    \centering
    \includegraphics[width=1\linewidth]{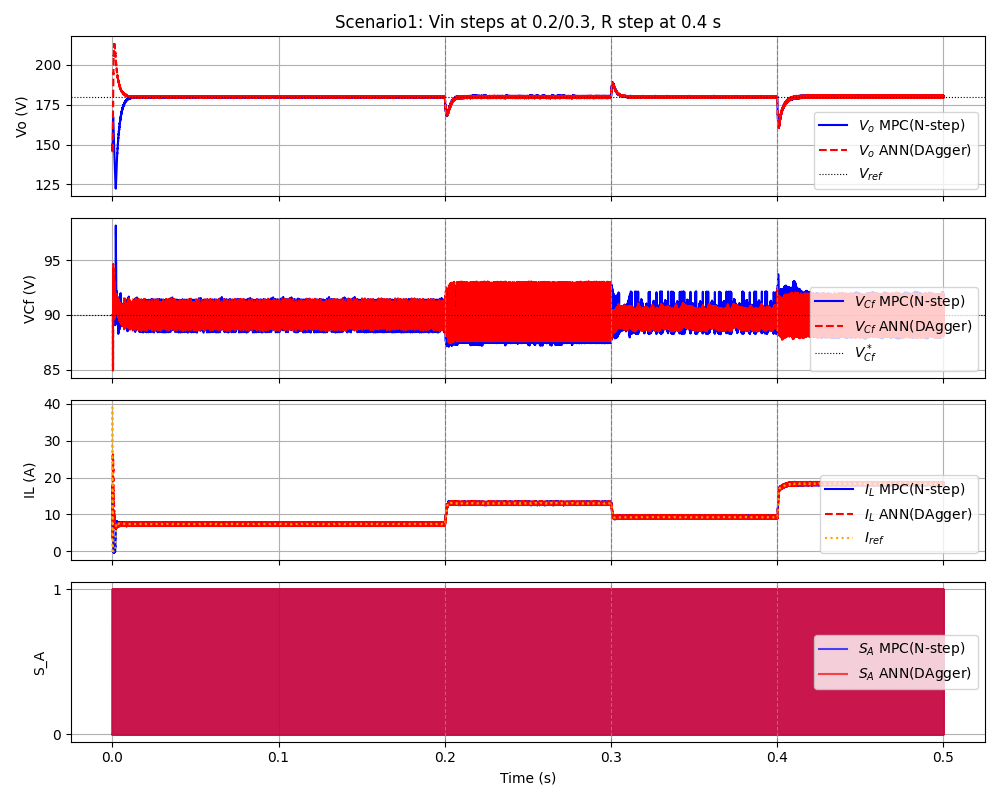}
    \caption{Closed-loop responses in Scenario~1 under nominal operating conditions, with input-voltage steps at 0.2~s and 0.3~s and a load-resistance step at 0.4~s. From top to bottom: output voltage \(v_o\), flying-capacitor voltage \(v_{Cf}\), inductor current \(i_L\), and switching signal \(S_A\). The distilled ANN closely matches the beam-search FCS-MPC expert.}
    \label{fig:scenario1}
\end{figure}

Scenario~1 considers nominal \(V_{\mathrm{in}}\) and load conditions with representative step disturbances. As shown in Fig.~\ref{fig:scenario1}, the closed-loop responses of the distilled ANN and the beam-search FCS-MPC expert are visually close. Quantitatively, Table~\ref{tab:key_results_all} shows that the ANN reduces \(\mathrm{MSE}_{v_o}\) from 14.13 to 6.22 and \(\mathrm{MSE}_{i_L}\) from 0.206 to 0.096, while preserving zero inductor-current violations. The ANN also lowers \(\mathrm{Overshoot}_{v_{Cf}}\) from 8.16~V to 4.65~V. The main trade-off is the output-voltage transient, where \(\mathrm{Overshoot}_{v_o}\) increases from 8.93~V to 33.39~V. Overall, Scenario~1 shows that the distilled policy reproduces the nominal closed-loop behavior of the long-horizon expert, with output-voltage overshoot as the main penalty.

\subsubsection{Scenario 2: Randomized Input Voltage and Load}
\label{sec:s2}

\begin{figure}
    \centering
    \includegraphics[width=1\linewidth]{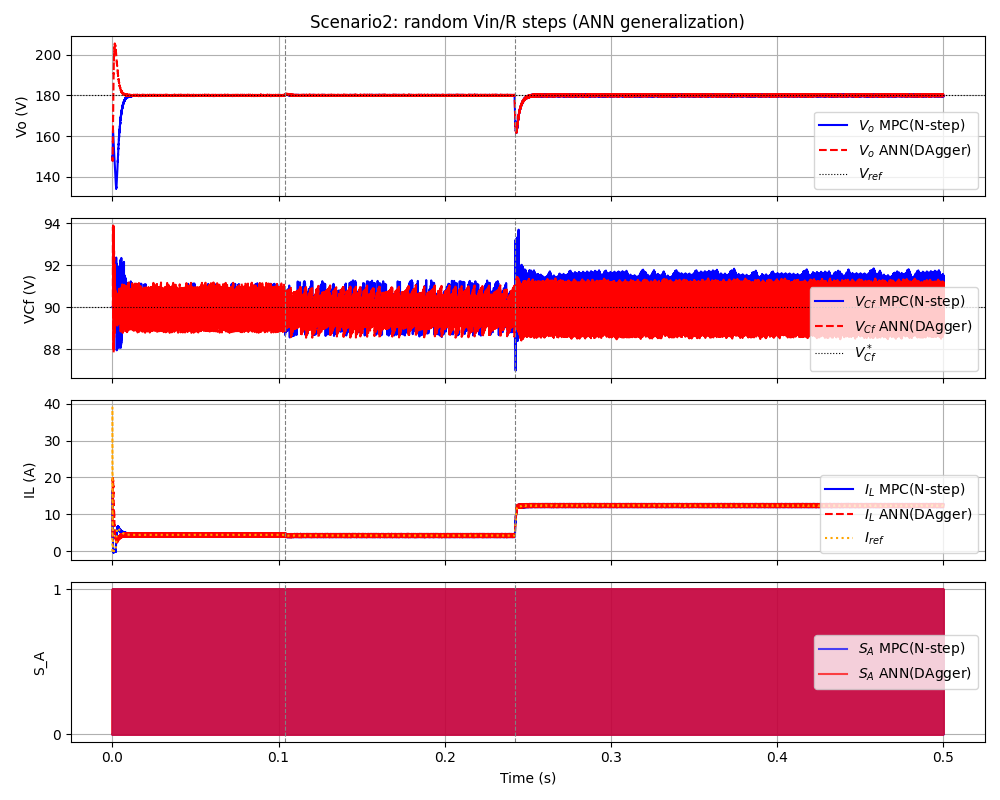}
    \caption{Closed-loop responses in Scenario~2 under randomized operating-point changes in input voltage and load resistance. From top to bottom: output voltage \(v_o\), flying-capacitor voltage \(v_{Cf}\), inductor current \(i_L\), and switching signal \(S_A\). The ANN policy remains stable and closely follows the FCS-MPC expert across varying operating conditions.}
    \label{fig:scenario2}
\end{figure}

Scenario~2 evaluates generalization under randomized step changes in \(V_{\mathrm{in}}\) and load within the domain-randomization ranges. As shown in Fig.~\ref{fig:scenario2}, the ANN remains stable across all operating intervals and continues to follow the expert closely at the waveform level. Table~\ref{tab:key_results_all} shows that the ANN reduces \(\mathrm{MSE}_{v_o}\) from 10.81 to 4.03 and \(\mathrm{MSE}_{i_L}\) from 0.175 to 0.076, while again maintaining zero inductor-current violations. The main discrepancy remains the transient output-voltage behavior, where \(\mathrm{Overshoot}_{v_o}\) increases from 0.69~V to 25.49~V. Thus, under operating-point randomization, the proposed policy preserves stable regulation and good current tracking, with output-voltage overshoot remaining the main trade-off.

\subsubsection{Scenario 3: Parameter Perturbations and Operating-Point Jumps}
\label{sec:s3}

\begin{figure}
    \centering
    \includegraphics[width=1\linewidth]{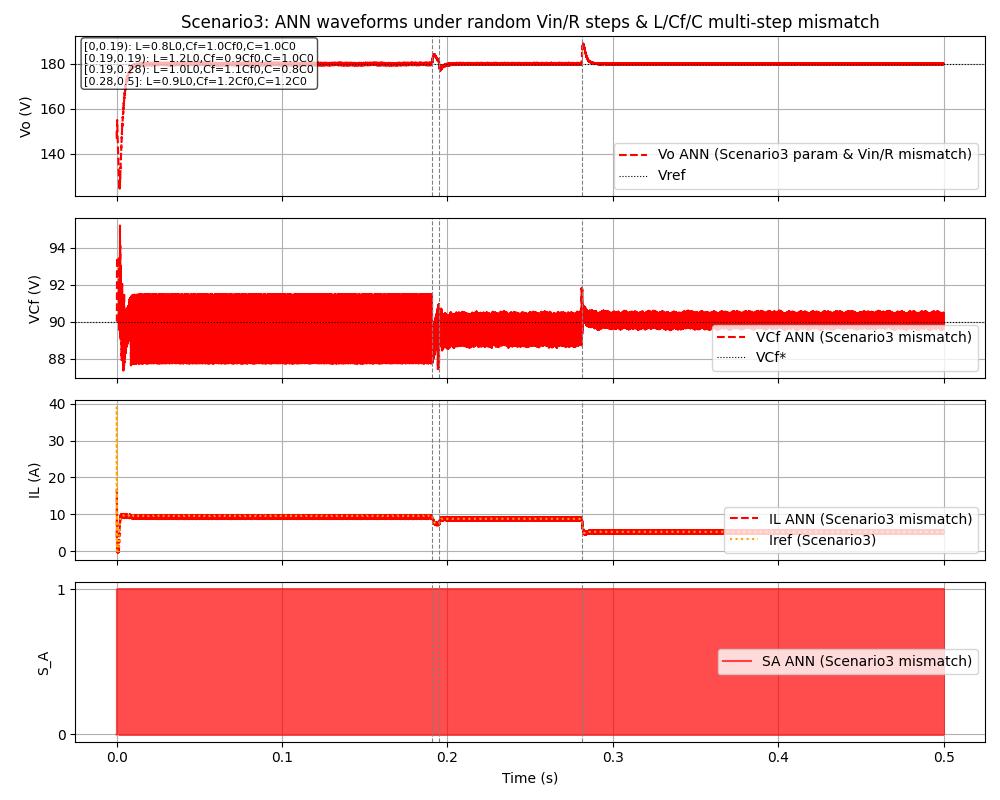}
    \caption{Representative ANN closed-loop responses in Scenario~3 under simultaneous operating-point randomization and passive-component perturbations in \(L\), \(C_f\), and \(C\). From top to bottom: output voltage \(v_o\), flying-capacitor voltage \(v_{Cf}\), inductor current \(i_L\), and switching signal \(S_A\). Stable regulation and capacitor balancing are preserved despite parametric mismatch.}
    \label{fig:scenario3}
\end{figure}

Scenario~3 further extends Scenario~2 by introducing passive-component perturbations in \((L, C_f, C)\) in addition to randomized \(V_{\mathrm{in}}\) and \(R\), making it the most demanding robustness case. As shown in Fig.~\ref{fig:scenario3}, the ANN still maintains stable closed-loop operation despite the combined operating-point shifts and plant mismatch. Table~\ref{tab:key_results_all} shows that the errors increase relative to Scenarios~1 and~2, with \(\mathrm{MSE}_{v_o}=33.94\), \(\mathrm{MSE}_{v_{Cf}}=5.170\), and \(\mathrm{Overshoot}_{v_o}=48.28~\mathrm{V}\). Nevertheless, all responses remain bounded, capacitor balancing is preserved, and no inductor-current violation occurs. These results support the robustness of the proposed training pipeline beyond nominal modeling assumptions.

\subsubsection{Training Summary, Inference Speed, and Objective Fidelity}

The ANN policy is trained using \(203{,}998\) MPC-labeled state--mode pairs generated under domain randomization and then refined by aggregating an additional \(50{,}000\) mismatch states collected with Disagreement-Based DAgger. After refinement, the classifier reaches a validation accuracy of \(0.9174\) and a test accuracy of \(0.9196\).

To quantify the computational savings, we measure per-step decision time for the \(N\)-step beam-search FCS-MPC expert and for the ANN policy on the same evaluation CPU. The ANN requires \(18.30~\mu\mathrm{s}\) per decision, whereas the expert requires \(342.26~\mu\mathrm{s}\), corresponding to an \(18.7\times\) speedup. Relative to the nominal control-update period of \(20~\mu\mathrm{s}\), the prototype ANN latency is slightly smaller. Since this timing is measured for a PyTorch prototype on the specific platform (Apple M3 Max CPU), it should be interpreted as a software-level runtime indicator rather than as a definitive guarantee of embedded real-time deployment. Nevertheless, the result confirms a substantial reduction in online decision cost.

To assess how closely the distilled policy reproduces the MPC objective, we compute the accumulated MPC stage cost \eqref{eq:mpc_cost} \emph{a posteriori} along the realized closed-loop trajectories and report the accumulated cost \(J_{\mathrm{sum}}\) and its per-step average \(J_{\mathrm{mean}}\) for Scenarios~1 and~2.
\begin{table}[t]
\centering
\caption{Objective-fidelity metrics under the same MPC cost \eqref{eq:mpc_cost}}
\label{tab:obj_fidelity}
\begin{tabular}{lcc}
\toprule
Scenario / Controller & \(J_{\mathrm{sum}}\) & \(J_{\mathrm{mean}}\) \\
\midrule
S1: MPC & \(1.0908\times10^{4}\) & 0.2182 \\
S1: ANN (DAgger) & \(5.3952\times10^{3}\) & 0.1079 \\
S2: MPC & \(9.0884\times10^{3}\) & 0.1818 \\
S2: ANN (DAgger) & \(4.1618\times10^{3}\) & 0.0832 \\
\bottomrule
\end{tabular}
\end{table}
Table~\ref{tab:obj_fidelity} shows that the ANN yields a lower realized accumulated cost than the beam-search expert in both Scenarios~1 and~2. We interpret this result cautiously. The expert is only an approximate solver because beam search with width \(K=15\) may prune switching sequences that would have achieved lower cumulative cost over the full horizon. In addition, disagreement-based DAgger retrains the student on learner-visited mismatch states, which can improve on-policy behavior in regions that were weakly represented in the original offline dataset. The ANN may also generate smoother switching sequences than the stepwise approximate expert. At the same time, the ANN still exhibits substantially larger output-voltage overshoot than the expert, so the lower realized cost should not be interpreted as uniformly better closed-loop control.

\subsection{Ablation Study}
Here, we report representative ablation results. To isolate the individual contributions of expert supervision, DR, and Disagreement-Based DAgger, we compare the four training configurations listed in Table~\ref{tab:ablation_cfg}. In particular, NO\_DR removes only the randomized offline data while retaining the same DAgger refinement, so that the effect of DR can be separated from the effect of on-policy correction. All configurations share the same ANN architecture, optimizer, \(20~\mu\mathrm{s}\) control-update period, and training schedule, and the Scenario~2 and Scenario~3 test trajectories are fixed across all configurations.

\begin{table}[t]
\centering
\caption{Ablation configurations and enabled components}
\label{tab:ablation_cfg}
\begin{tabular}{lccc}
\toprule
Config & Expert Labels & DR & DAgger \\
\midrule
FULL       & \checkmark & \checkmark & \checkmark \\
NO\_DAGGER & \checkmark & \checkmark & \(\times\) \\
NO\_DR     & \checkmark & \(\times\) & \checkmark \\
NO\_EXPERT & \(\times\) & N/A & N/A \\
\bottomrule
\end{tabular}
\end{table}

\begin{table}[t]
\centering
\caption{Representative ablation results under Scenarios~1--3 (lower is better for all metrics shown)}
\label{tab:ablation_results}
\scriptsize
\setlength{\tabcolsep}{3pt}
\renewcommand{\arraystretch}{1.05}
\resizebox{\columnwidth}{!}{%
\begin{tabular}{@{}llrrrr@{}}
\toprule
Scenario & Metric & FULL & NO\_DAGGER & NO\_DR & NO\_EXPERT \\
\midrule
S1 & \(\mathrm{MSE}_{v_o}\) & \textbf{13.9237} & 14.0757 & 14.1253 & 15666.8631 \\
   & \(\mathrm{MSE}_{i_L}\) & \textbf{0.2118} & 0.2885 & 0.2669 & 2461.4469 \\
   & \(\mathrm{Overshoot}_{v_o}\) & \textbf{7.0795} & 7.6754 & 7.2936 & 738.6373 \\
   & \(N_{i_L,\mathrm{viol}}\) & \textbf{0} & \textbf{0} & \textbf{0} & 1939 \\
\midrule
S2 & \(\mathrm{MSE}_{v_o}\) & 13.2922 & 13.9687 & \textbf{13.1520} & 219519.8731 \\
   & \(\mathrm{MSE}_{i_L}\) & \textbf{0.1954} & 1.0585 & 0.5263 & 12721.0718 \\
   & \(\mathrm{Overshoot}_{v_o}\) & \textbf{10.6224} & 12.4090 & 15.6412 & 720.7855 \\
   & \(N_{i_L,\mathrm{viol}}\) & \textbf{0} & \textbf{0} & \textbf{0} & 40333 \\
\midrule
S3 & \(\mathrm{MSE}_{v_o}\) & \textbf{8.5146} & 8.6940 & 16.0551 & 380609.4699 \\
   & \(\mathrm{MSE}_{i_L}\) & 0.2935 & \textbf{0.2815} & 23.0244 & 18813.3168 \\
   & \(\mathrm{Overshoot}_{v_o}\) & \textbf{5.3896} & 5.9888 & 10.9562 & 877.9873 \\
   & \(N_{i_L,\mathrm{viol}}\) & \textbf{0} & \textbf{0} & \textbf{0} & 42861 \\
\bottomrule
\end{tabular}%
}
\end{table}

Rather than imposing a strict total ordering across all scenarios and metrics, Table~\ref{tab:ablation_results} supports three robust conclusions. First, expert supervision is indispensable. Second, DR is the main source of robustness beyond nominal conditions. Third, Disagreement-Based DAgger provides additional gains mainly in on-policy current-tracking and transient behavior.

The importance of expert supervision is most clearly seen from the NO\_EXPERT configuration. This setting fails to produce a viable closed-loop policy in all three scenarios, with errors increasing by orders of magnitude and thousands of current-limit violations. For example, \(\mathrm{MSE}_{v_o}\) rises to \(15666.8631\), \(219519.8731\), and \(380609.4699\) in Scenarios~1, 2, and~3, respectively. These results confirm that MPC-derived expert labels are essential for learning a stabilizing switching policy under the present network architecture and training setup.

The role of DR becomes clear by comparing FULL and NO\_DR. Under nominal conditions (Scenario~1), the three expert-supervised configurations remain close to one another, indicating that nominal-data training is sufficient when the training and test distributions are well matched. Under operating-point randomization (Scenario~2), NO\_DR remains stable, but its current-tracking and transient metrics degrade relative to FULL. Although NO\_DR attains a slightly smaller \(\mathrm{MSE}_{v_o}\) than FULL in Scenario~2, it exhibits worse \(\mathrm{MSE}_{i_L}\) and larger output-voltage overshoot. The effect of DR becomes much more pronounced in Scenario~3, where operating-point variation is combined with passive-component perturbations: relative to FULL, NO\_DR increases \(\mathrm{MSE}_{i_L}\) from 0.2935 to 23.0244, \(\mathrm{MSE}_{v_o}\) from 8.5146 to 16.0551, and \(\mathrm{Overshoot}_{v_o}\) from 5.3896 to 10.9562. These results indicate that DR is the primary mechanism enabling robustness to joint operating-point shifts and parameter mismatch.

The contribution of Disagreement-Based DAgger is isolated by comparing FULL with NO\_DAGGER. In Scenario~1, the two are close, although FULL still improves \(\mathrm{MSE}_{i_L}\) and slightly reduces output-voltage overshoot. The clearest gains appear in Scenario~2, where FULL reduces \(\mathrm{MSE}_{i_L}\) from 1.0585 to 0.1954 and \(\mathrm{Overshoot}_{v_o}\) from 12.4090 to 10.6224. In Scenario~3, the difference is more nuanced: NO\_DAGGER slightly improves \(\mathrm{MSE}_{i_L}\), but FULL achieves lower \(\mathrm{MSE}_{v_o}\) and lower output-voltage overshoot. This suggests that under the strongest perturbations, DAgger mainly improves transient quality and suppresses extreme on-policy deviations, even when some average-error metrics are already comparable.

Overall, the ablation study shows that all expert-supervised models perform similarly under nominal conditions, DR is the dominant factor that preserves robustness under randomized operating conditions and parameter perturbations, and Disagreement-Based DAgger yields additional benefits once the learner visits states that are weakly represented in the original offline dataset. Sensitivity experiments examining the DAgger mismatch-sample budget and DR intensity are reported in Appendix~\ref{sec:sensitivity}.
\subsection{Transfer Learning Experiments}
\label{sec:transfer}

\begin{figure}
    \centering
    \includegraphics[width=1\linewidth]{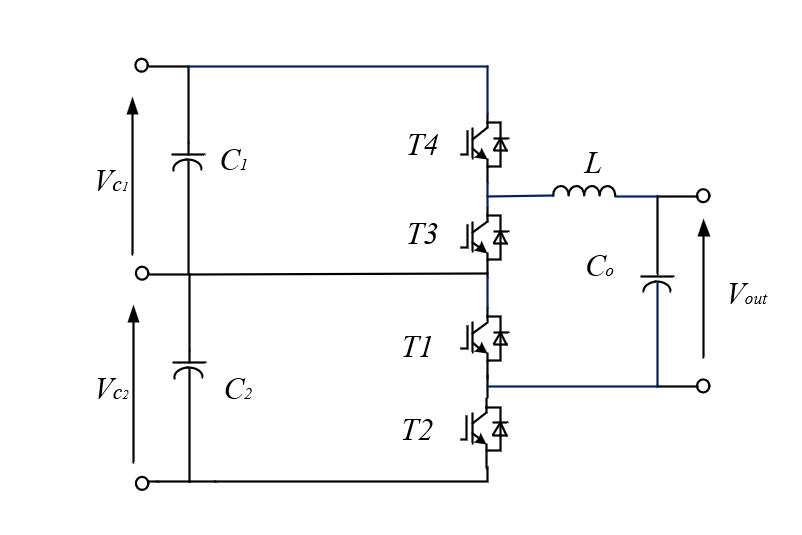}
    \caption{Topology of NPC-Buck Converter}
    \label{fig:npc_topology}
\end{figure}

A natural question is whether the neural features learned for one converter topology can be reused for a related but distinct topology, thereby reducing the data and training effort required for the new target system.
The FC-TLBC and the NPC-type three-level buck converter (Buck-3L) share several structural properties that make cross-topology transfer plausible.
First, both are three-level converter topologies whose switching behavior can be described by the same number of discrete modes \((|\mathcal{U}|=4)\).
Second, the state vectors in both cases consist of an inductor current, an internal capacitor voltage, and an output voltage, so the six-dimensional input feature space \(\mathbf{z}_k\) defined in~\eqref{eq:feature_vector} has the same physical interpretation.
Third, the control objectives---current tracking and internal capacitor-voltage balancing subject to mode-feasibility constraints---are analogous, differing mainly in the sign conventions and mode-coefficient values of the state-space matrices.
These commonalities suggest that the hidden-layer weights trained on FC-TLBC data already encode useful nonlinear decision boundaries that are transferable to Buck-3L with only output-layer adaptation.

To evaluate this hypothesis, we follow the protocol in Algorithm~\ref{alg:transfer}, which separates source pre-training on FC-TLBC, target training from scratch on Buck-3L, and transfer initialization/fine-tuning on the same Buck-3L target dataset.

\begin{algorithm}[t]
\caption{Transfer-Learning Evaluation Protocol}
\label{alg:transfer}
\begin{algorithmic}[1]
\State \textbf{Stage 1: Source pre-training on FC-TLBC}
\State Generate the FC-TLBC expert dataset \(\mathcal{D}_{\mathrm{FC}}\)
\State Train a source ANN policy \(\pi_{\theta_{\mathrm{src}}}\) on \(\mathcal{D}_{\mathrm{FC}}\)
\State Evaluate source-domain accuracy to confirm convergence
\Statex
\State \textbf{Stage 2: Buck-3L training from scratch}
\State Generate the Buck-3L expert dataset \(\mathcal{D}_{\mathrm{Buck}}\)
\State Initialize \(\pi_{\theta_{\mathrm{scratch}}}\) with random weights
\State Train \(\pi_{\theta_{\mathrm{scratch}}}\) on \(\mathcal{D}_{\mathrm{Buck}}\)
\State Evaluate closed-loop metrics under the Buck-3L test scenarios
\Statex
\State \textbf{Stage 3: Buck-3L transfer learning}
\State Initialize \(\pi_{\theta_{\mathrm{trans}}}\) from \(\pi_{\theta_{\mathrm{src}}}\)
\State Re-initialize only the output layer for the Buck-3L action set
\State Fine-tune \(\pi_{\theta_{\mathrm{trans}}}\) on \(\mathcal{D}_{\mathrm{Buck}}\)
\State Compare MPC, Scratch, and Transfer under the same Buck-3L scenarios
\end{algorithmic}
\end{algorithm}

We consider three controllers:
\begin{enumerate}
  \item \textbf{MPC:} FCS-MPC tailored for Buck-3L and used as the reference controller.
  \item \textbf{Scratch:} Buck-3L controller trained from random initialization using 4053 MPC-labeled samples and 40 epochs.
  \item \textbf{Transfer:} Initialize the Buck-3L network with hidden-layer weights from an FC-TLBC source model trained specifically for this transfer experiment on 8203 source samples for 60 epochs; re-initialize only the output layer; then fine-tune on the same 4053 Buck samples for 40 epochs.
\end{enumerate}

For this transfer-learning study, the FC-TLBC source model described above achieves a test accuracy of about 0.86 on its source-domain split.
On Buck-3L, the Scratch model reaches a test accuracy of 0.80--0.83, while the Transfer model reaches approximately 0.94, indicating that the source-domain features improve action classification with the same amount of target data.

Closed-loop performance is evaluated under two step-load scenarios S1 and S2, with reference voltage \(v_o^\star=80\)~V, input voltage around 120~V, and load resistance stepping from \(20\,\Omega\) to \(10\,\Omega\):
\begin{itemize}
  \item In S1 (moderate disturbance), MPC and Transfer responses almost overlap, with peak overshoot \(\approx 0.2\)~V (\(0.22\%\)), while Scratch exhibits noticeable oscillation and larger \(\mathrm{MSE}_{v_o}\) (7.74 vs. 3.95 for MPC and 3.71 for Transfer).
  \item In S2 (strong disturbance), Scratch yields severe over-voltage (up to about 120~V, \(88.9\%\) overshoot) and slow recovery, with \(\mathrm{MSE}_{v_o}=336.9\). Transfer maintains \(\mathrm{MSE}_{v_o}=3.01\) and overshoot \(\approx 3.24\)~V (\(4.05\%\)), close to MPC's \(\mathrm{MSE}_{v_o}=2.33\).
\end{itemize}

Average efficiency \(\mathrm{Eff}_{\mathrm{avg}}\) and average output power \(P_{\mathrm{out,avg}}\) are similar across MPC, Scratch, and Transfer, indicating that improved tracking does not come at the cost of energy efficiency.

\begin{figure}
    \centering
    \includegraphics[width=1\linewidth]{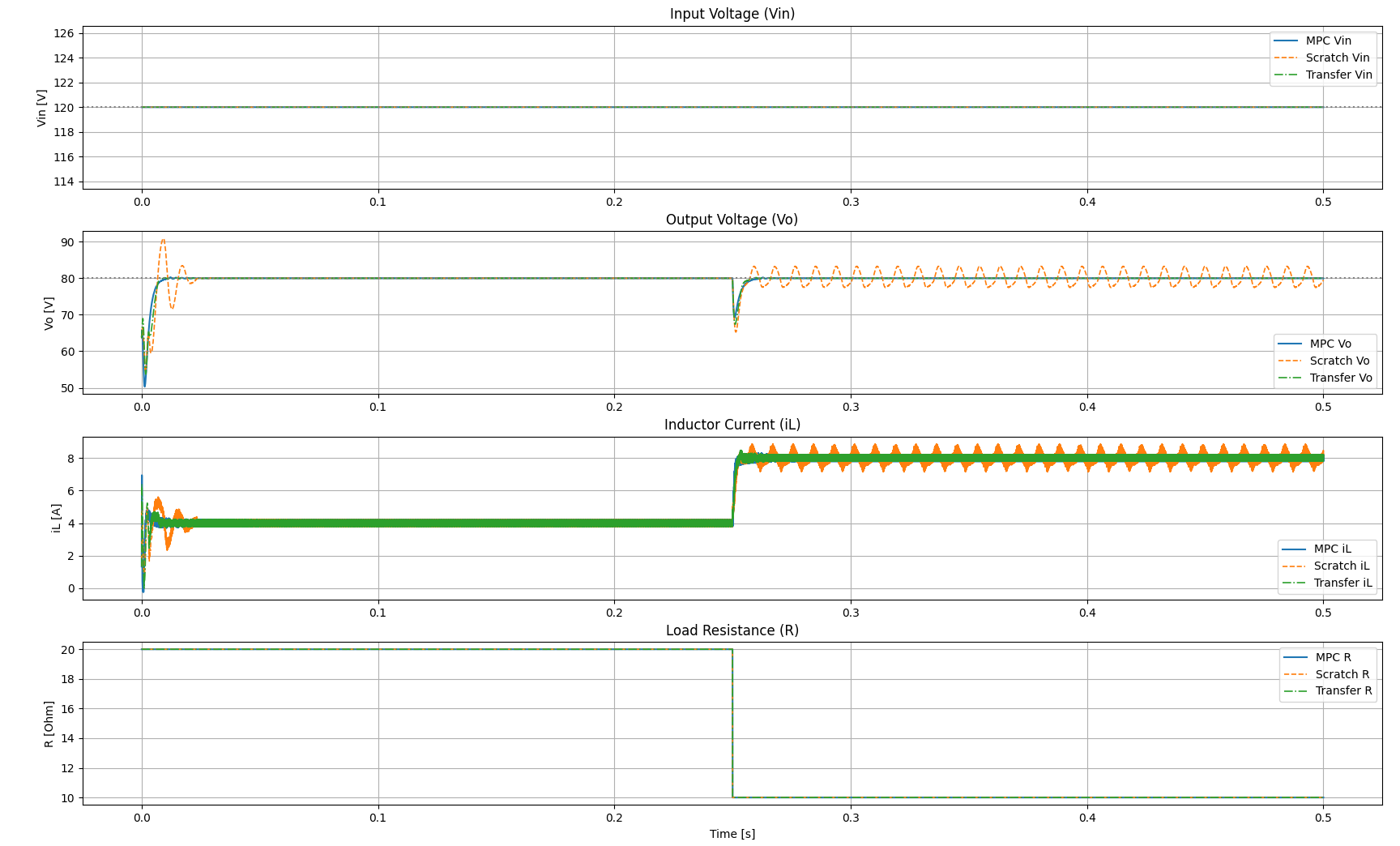}
    \caption{Closed-loop comparison in Transfer Learning Scenario~1 for the NPC-type three-level buck converter under a moderate load step. From top to bottom: input voltage, output voltage, inductor current, and load resistance. The transferred policy closely follows the Buck-3L MPC reference and improves over training from scratch.}
    \label{fig:transfer1}
\end{figure}

\begin{figure}
    \centering
    \includegraphics[width=1\linewidth]{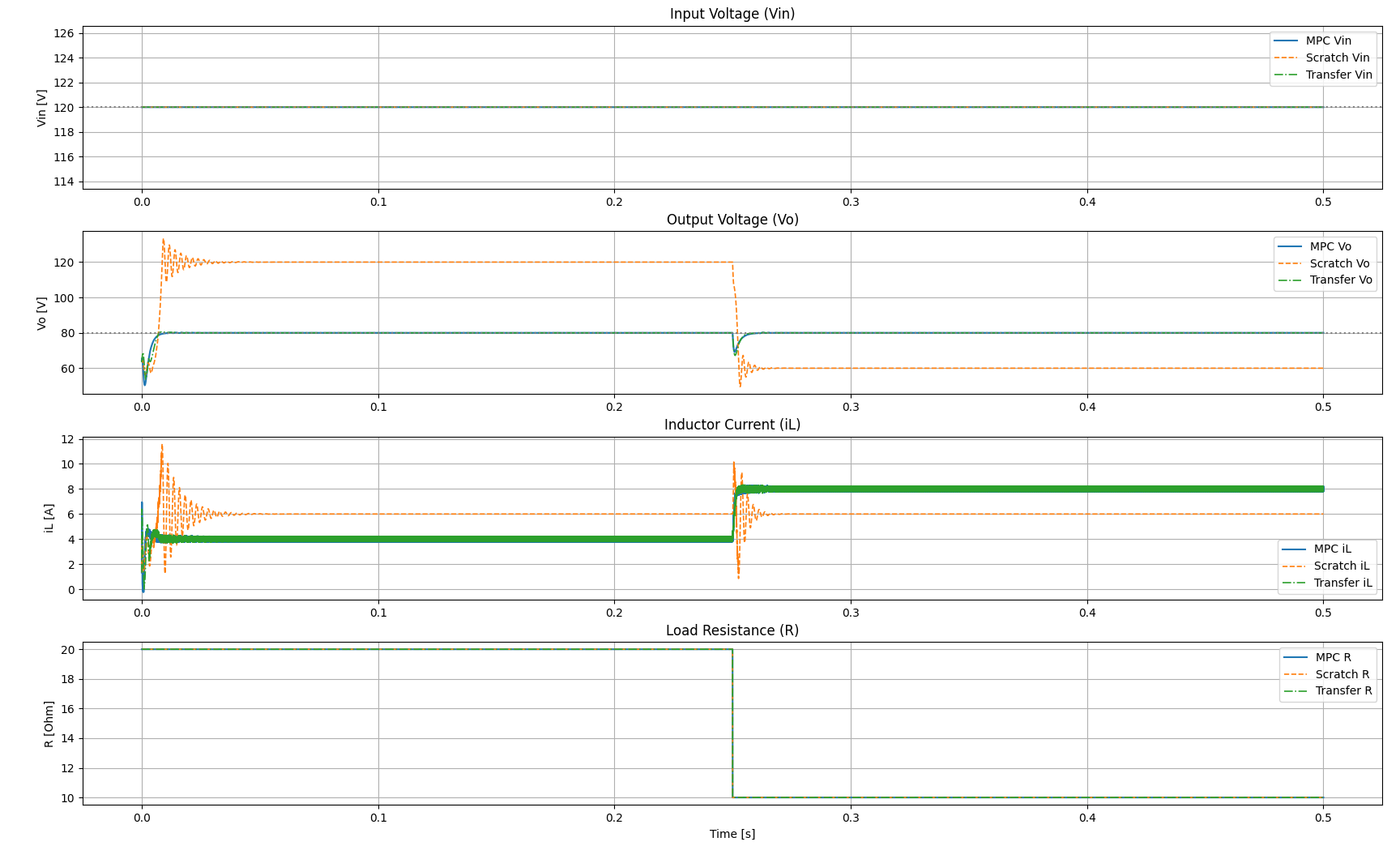}
    \caption{Closed-loop comparison in Transfer Learning Scenario~2 for the NPC-type three-level buck converter under a larger disturbance. From top to bottom: input voltage, output voltage, inductor current, and load resistance. The transferred policy remains close to the MPC reference, whereas the scratch-trained policy exhibits larger overshoot and slower recovery.}
    \label{fig:transfer2}
\end{figure}

These results confirm that:
\begin{itemize}
  \item features learned on FC-TLBC are reusable on Buck-3L,
  \item transfer learning improves Buck-3L performance with the same data budget, and
  \item cross-topology generalization is feasible within the proposed MPC-to-ANN framework.
\end{itemize}

\begin{table}[t]
\centering
\caption{Controllers Compared in Transfer Learning Experiments}
\label{tab:transfer_cfg}
\begin{tabular}{l l}
\toprule
Controller & Description \\
\midrule
MPC & FCS-MPC expert tailored for Buck-3L \\
Scratch & Buck-3L ANN trained from random initialization \\
Transfer & Buck-3L ANN initialized from FC-TLBC source model \\
\bottomrule
\end{tabular}
\end{table}
\section{Conclusion}
This paper presented a practical MPC-to-neural distillation framework for FC-TLBCs, where a compact feedforward switching policy is learned from a long-horizon beam-search FCS-MPC expert. By combining domain-randomized expert demonstrations with disagreement-based DAgger refinement, the proposed method reduces the online computational burden while improving robustness to operating-point variation and passive-component mismatch.

Simulation results showed that the distilled controller preserves stable output-voltage regulation and flying-capacitor balancing under nominal conditions, randomized operating points, and parameter perturbations. On the evaluation CPU, the per-decision computation time was reduced. 
The main limitation may be that the ANN exhibits larger output-voltage overshoot than the MPC expert in Scenarios~1 and~2. The ablation study further showed that expert supervision is essential, domain randomization is the main driver of robustness, and disagreement-based DAgger yields additional gains in on-policy transient and current-tracking behavior.

The transfer-learning results suggest that representations learned on FC-TLBC can be reused for a related three-level buck topology, improving data efficiency relative to training from scratch. Future work will focus on embedded and experimental validation and on extending the training pipeline to account for nonideal effects such as dead time, switching losses, and measurement noise. Overall, the results indicate that neural distillation is a practical route for bringing long-horizon predictive control closer to real-time use in multilevel power converters.

 
\bibliographystyle{ieeetr}
\nocite{*}
\bibliography{References_revised_full} 
 
\appendices
\section{Metrics Used in the Experiments}
\label{sec:metric_defs}

All trajectory-based metrics are evaluated over a closed-loop rollout of
\(N_{\mathrm{sim}}\) samples with control-update period \(T_s\).
We define
\[
T_{\mathrm{total}} = N_{\mathrm{sim}} T_s,
\]
and let \(t_k\) denote the physical time associated with sample \(k\).
The voltage references are
\[
V_{\mathrm{ref}} = v_o^\star, \qquad
V_{Cf,\mathrm{ref}} = \frac{v_o^\star}{2},
\]
and \(i_{\mathrm{ref},k}\) is generated by the outer voltage controller.

The reported tracking and transient metrics are defined as follows:
\begin{align}
\mathrm{MSE}_{v_o}
&=
\frac{1}{N_{\mathrm{sim}}}
\sum_{k=1}^{N_{\mathrm{sim}}}
\bigl(v_{o,k}-V_{\mathrm{ref}}\bigr)^2, \\
\mathrm{MSE}_{v_{Cf}}
&=
\frac{1}{N_{\mathrm{sim}}}
\sum_{k=1}^{N_{\mathrm{sim}}}
\bigl(v_{Cf,k}-V_{Cf,\mathrm{ref}}\bigr)^2, \\
\mathrm{MSE}_{i_L}
&=
\frac{1}{N_{\mathrm{sim}}}
\sum_{k=1}^{N_{\mathrm{sim}}}
\bigl(i_{L,k}-i_{\mathrm{ref},k}\bigr)^2.
\end{align}

We also report the signed final-sample steady-state error:
\begin{align}
\mathrm{SSE}_{v_o}
&=
v_{o,N_{\mathrm{sim}}}-V_{\mathrm{ref}}, \\
\mathrm{SSE}_{v_{Cf}}
&=
v_{Cf,N_{\mathrm{sim}}}-V_{Cf,\mathrm{ref}}.
\end{align}

The peak overshoot and its percentage form are defined by
\begin{align}
\mathrm{Overshoot}_{v_o}
&=
\max_{1 \le k \le N_{\mathrm{sim}}} v_{o,k} - V_{\mathrm{ref}}, \\
\mathrm{Overshoot}_{v_{Cf}}
&=
\max_{1 \le k \le N_{\mathrm{sim}}} v_{Cf,k} - V_{Cf,\mathrm{ref}},
\end{align}
and
\begin{align}
M_{p,v_o}(\%)
&=
100\frac{\mathrm{Overshoot}_{v_o}}{V_{\mathrm{ref}}}, \\
M_{p,v_{Cf}}(\%)
&=
100\frac{\mathrm{Overshoot}_{v_{Cf}}}{V_{Cf,\mathrm{ref}}}.
\end{align}

The settling times are computed using a \(\pm 2\%\) band:
\begin{align*}
T_{\mathrm{set},v_o}
&=
\max\left\{
t_k \,\middle|\,
v_{o,k}\notin [0.98V_{\mathrm{ref}},\,1.02V_{\mathrm{ref}}]
\right\}, \\
T_{\mathrm{set},v_{Cf}}
&=
\max\left\{
t_k \,\middle|\,
v_{Cf,k}\notin [0.98V_{Cf,\mathrm{ref}},\,1.02V_{Cf,\mathrm{ref}}]
\right\}.
\end{align*}
For the multi-step scenarios considered here,
\(T_{\mathrm{set},v_o}\) and \(T_{\mathrm{set},v_{Cf}}\)
should therefore be interpreted as the last-exit time from the
\(\pm 2\%\) band over the entire rollout.

The steady-state ripple is evaluated as the standard deviation after
\(t \ge 0.4~\mathrm{s}\):
\begin{align*}
\mathrm{Ripple}_{v_o}
&=
\mathrm{std}\!\left(
\{\,v_{o,k} \mid t_k \ge 0.4~\mathrm{s}\,\}
\right), \\
\mathrm{Ripple}_{v_{Cf}}
&=
\mathrm{std}\!\left(
\{\,v_{Cf,k} \mid t_k \ge 0.4~\mathrm{s}\,\}
\right).
\end{align*}

The over-voltage and sag penalties are defined as
\begin{align}
\mathrm{Penalty}_{\mathrm{over}}
&=
\frac{T_s}{V_{\mathrm{ref}}}
\sum_{k=1}^{N_{\mathrm{sim}}}
\max\!\left(v_{o,k}-1.05V_{\mathrm{ref}},\,0\right), \\
\mathrm{Penalty}_{\mathrm{sag}}
&=
\frac{T_s}{V_{\mathrm{ref}}}
\sum_{k=1}^{N_{\mathrm{sim}}}
\max\!\left(0.95V_{\mathrm{ref}}-v_{o,k},\,0\right).
\end{align}

The inductor-current violation count is
\begin{equation}
N_{i_L,\mathrm{viol}}
=
\sum_{k=1}^{N_{\mathrm{sim}}}
\mathbf{1}\!\left(
i_{L,k}\notin \mathcal{I}_{\mathrm{safe}}
\right),
\end{equation}
where \(\mathcal{I}_{\mathrm{safe}}\) denotes the hard current-limit interval
used in the controller design and simulator.

The switching statistics are defined by
\begin{align}
s_k &= [S_{A,k},\,S_{B,k}], \\
\mathrm{SwitchCount}
&=
\sum_{k=2}^{N_{\mathrm{sim}}}
\mathbf{1}\!\left(s_k \neq s_{k-1}\right), \\
\mathrm{SwitchFreq}
&=
\frac{\mathrm{SwitchCount}}{T_{\mathrm{total}}}, \\
N_{S_A}
&=
\sum_{k=2}^{N_{\mathrm{sim}}}
\mathbf{1}\!\left(S_{A,k}\neq S_{A,k-1}\right), \\
N_{S_B}
&=
\sum_{k=2}^{N_{\mathrm{sim}}}
\mathbf{1}\!\left(S_{B,k}\neq S_{B,k-1}\right), \\
N_{\mathrm{trans,total}}
&=
N_{S_A}+N_{S_B}.
\end{align}

If the energy-related quantities are retained, they are computed as
\begin{align}
E_{\mathrm{in}}
&=
T_s
\sum_{k=1}^{N_{\mathrm{sim}}}
V_{\mathrm{in},k}\,i_{L,k}, \\
E_{\mathrm{out}}
&=
T_s
\sum_{k=1}^{N_{\mathrm{sim}}}
v_{o,k}\,i_{o,k}, \\
P_{\mathrm{out,avg}}
&=
\frac{E_{\mathrm{out}}}{T_{\mathrm{total}}}, \\
\mathrm{Eff}_{\mathrm{avg}}
&=
\frac{E_{\mathrm{out}}}{E_{\mathrm{in}}}.
\end{align}

\section{Sensitivity Experiments}
\label{sec:sensitivity}

This appendix evaluates how sensitive the proposed learning pipeline is to
two key design choices: (i) the Disagreement-Based DAgger mismatch-sample
budget \(N_{\mathrm{Dag}}\) and (ii) the strength of domain randomization (DR)
used to generate the offline expert dataset. We focus on the eight
highest-variance metrics for each scenario, as these are the most informative
about what actually changes when \(N_{\mathrm{Dag}}\) or the DR intensity is varied.
 
\subsection{Disagreement-Based DAgger Sample Size Sensitivity}
 
Disagreement-Based DAgger's effect depends on the number of mismatch samples $N_{\text{Dag}}$. We evaluate $N_{\text{Dag}}\in\{0, 500,1000,2000,4000,8000,12000\}$, starting from the same DR-pretrained model. For each setting, we collect up to $N_{\text{Dag}}$ mismatch states in closed loop, retrain the network, and then evaluate on Scenarios 2 and 3.
 
\begin{figure}
    \centering
    \includegraphics[width=1\linewidth]{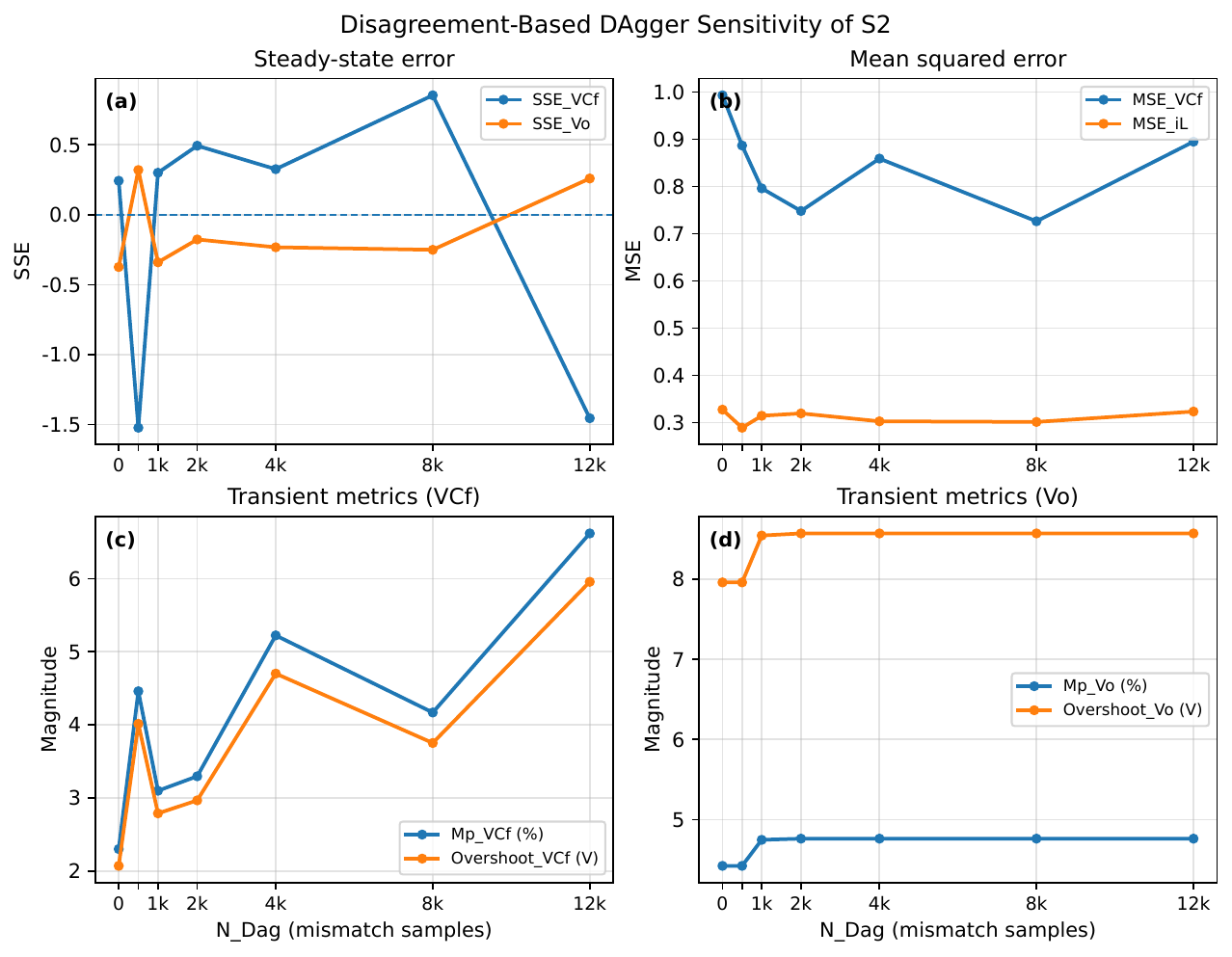}
    \caption{Disagreement-Based DAgger Sensitivity of S2}
    \label{fig:dagger_s2}
\end{figure}
\begin{figure}
    \centering
    \includegraphics[width=1\linewidth]{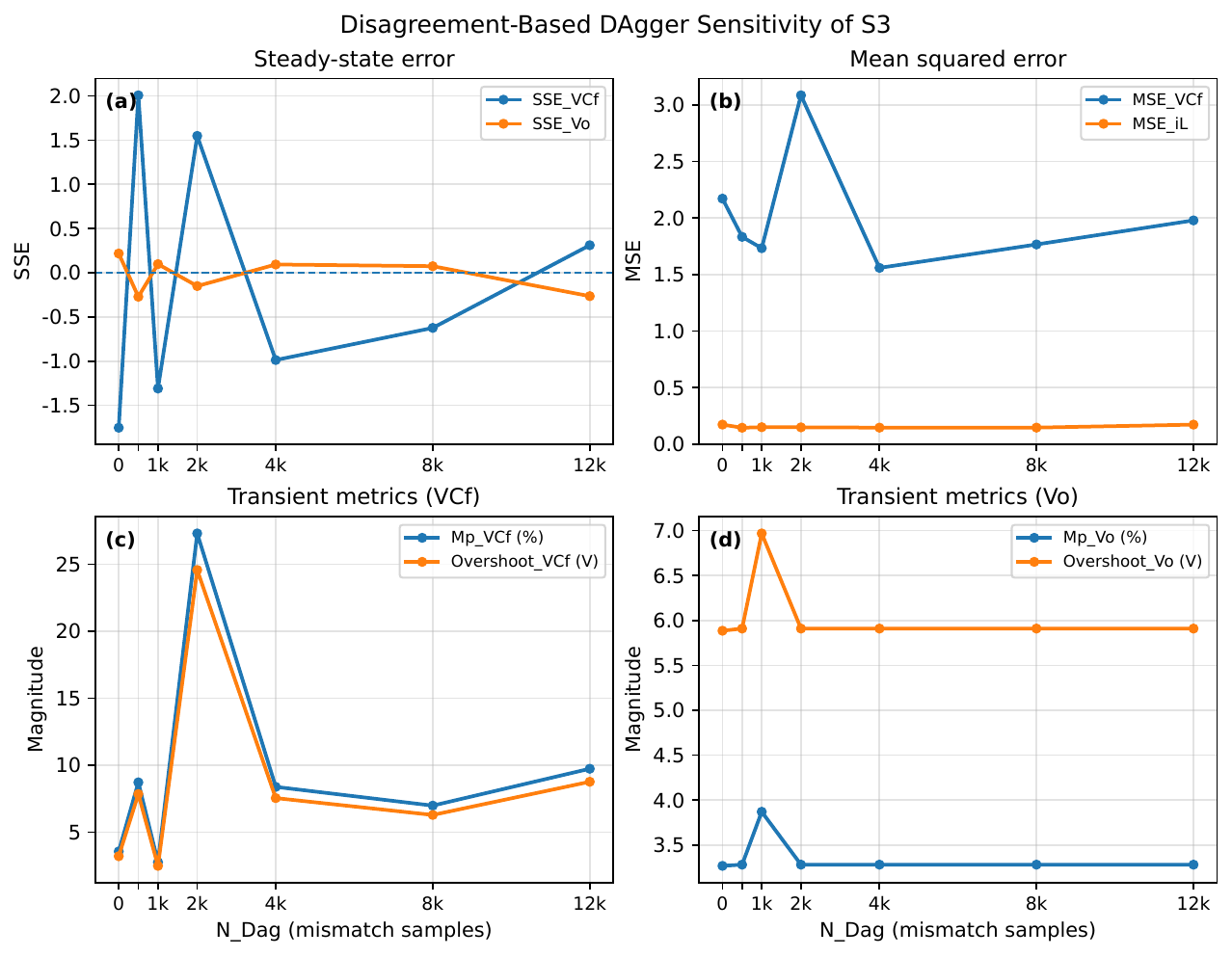}
    \caption{Disagreement-Based DAgger Sensitivity of S3}
    \label{fig:dagger_s3}
\end{figure}
 
The key observations are:
\begin{itemize}
  \item \textbf{Rapid changes in transient metrics with small budgets:}
  The most visibly moving curves are the peak/overshoot-related terms
  ($Mp_{Vcf,\%}$ and $Overshoot_{Vcf}$ in particular), indicating that adding a small number of disagreement samples mainly corrects
  \emph{switching-boundary and transient decisions}, reducing voltage spikes more than it changes steady tracking.
  \item \textbf{A practical stability region (few thousand samples):}
  For intermediate budgets ($N_{\text{Dag}}\approx 1000$--$8000$), the majority of the plotted metrics settle into a relatively stable range.
  \item \textbf{Non-monotonic behavior at very large budgets:}
  At $N_{\text{Dag}}=12000$, several transient-dominant metrics can rise again, consistent with mismatch states being over-represented near switching boundaries and the beam-search expert providing less consistent labels in rarely visited states.
\end{itemize}
 
This suggests that Disagreement-Based DAgger is highly sample-efficient: a few thousand additional expert queries are sufficient to obtain most of the improvement, especially in peak/overshoot behavior.
 
\subsection{Domain Randomization Intensity Sensitivity}
 
To evaluate DR intensity, we scale the randomization range as $r\in\{10\%,30\%,50\%,80\%,100\%\}$ relative to the full range used in the main experiments. For each $r$, we regenerate the DR dataset, retrain the ANN for 40 epochs, and evaluate on the fixed Scenario 2 and Scenario 3 test sets.

\begin{figure}
    \centering
    \includegraphics[width=1\linewidth]{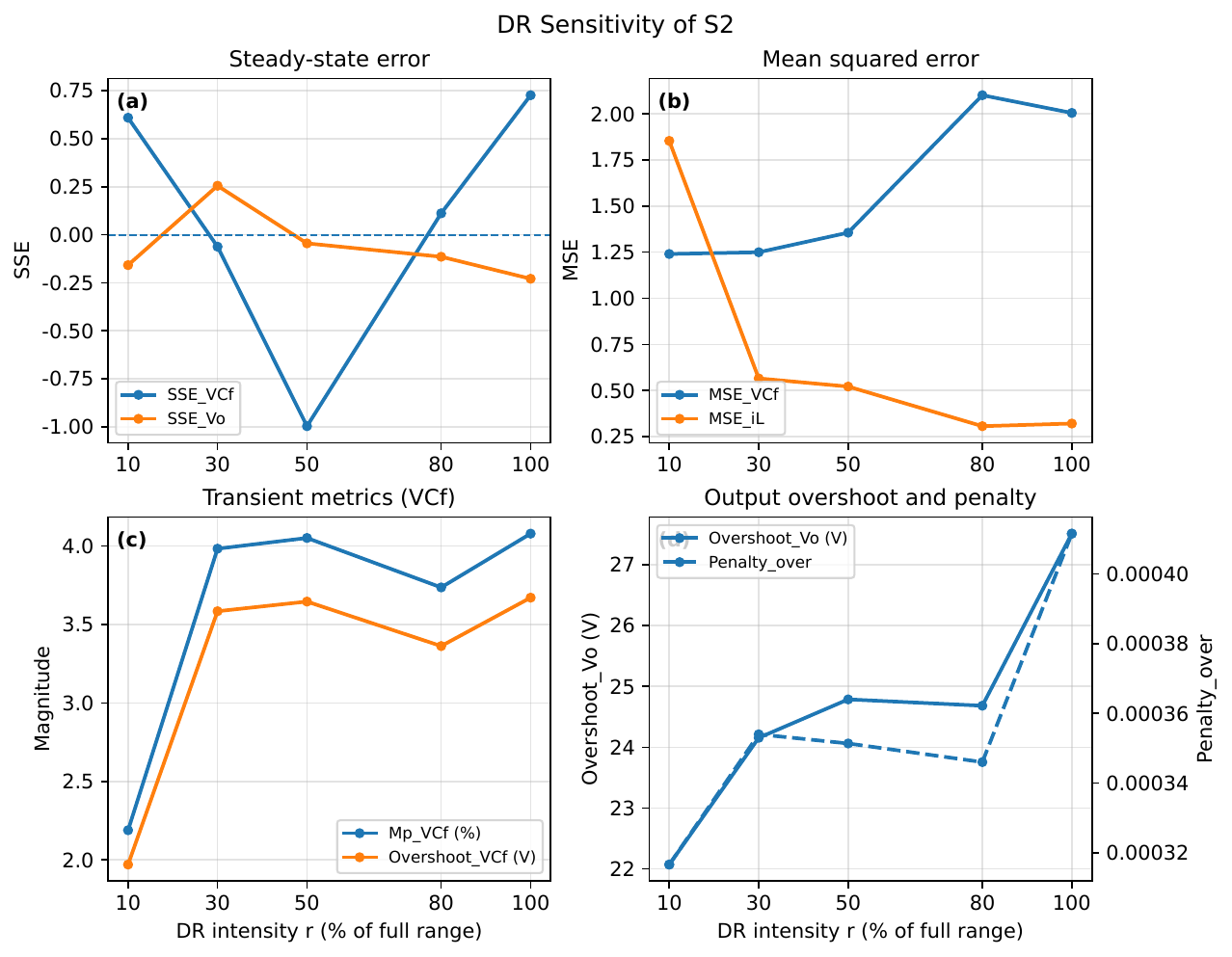}
    \caption{DR Sensitivity of S2}
    \label{fig:dr_s2}
\end{figure}
 
\begin{figure}
    \centering
    \includegraphics[width=1\linewidth]{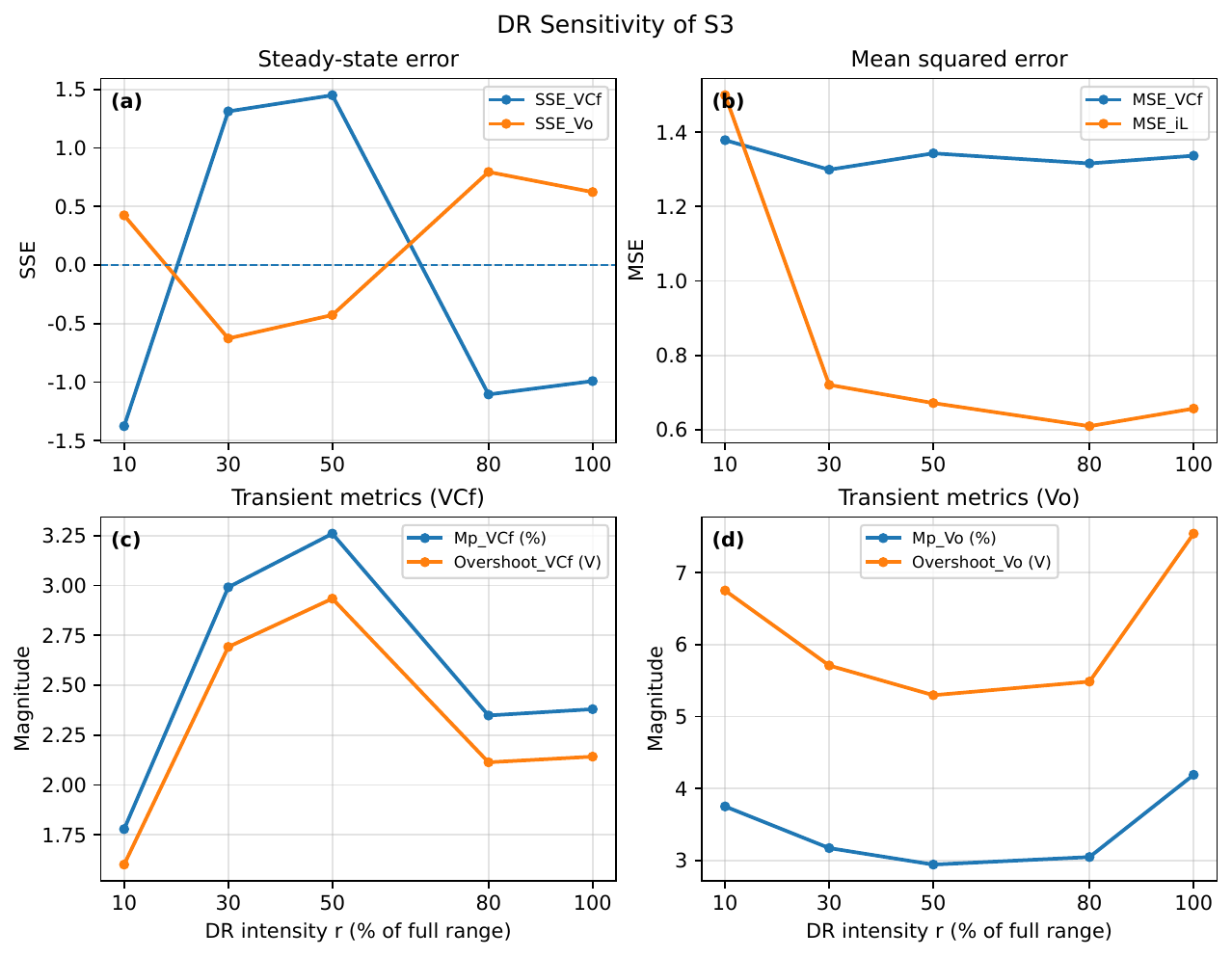}
    \caption{DR Sensitivity of S3}
    \label{fig:dr_s3}
\end{figure}
 
The results show:
\begin{itemize}
  \item \textbf{Under-coverage (10\% DR):}
  Insufficient randomization leads to poorer robustness, most evident in dynamic-tracking metrics ($MSE_{i_L}$ and $SSE$ terms).
  \item \textbf{Intermediate ranges (30\%--50\%):}
  The best trade-off is achieved at intermediate DR, keeping both average errors and transient measures in a balanced regime.
  \item \textbf{Very strong DR (80\%--100\%):}
  Increasing DR further does not necessarily improve the averages and can \emph{worsen transient behavior}, as the approximation task becomes harder.
  \item Average efficiency $Eff_{\text{avg}}$ is almost invariant across DR ranges, indicating that DR mainly affects dynamic tracking and not steady-state power conversion quality.
\end{itemize}
Overall, DR exhibits an ``intermediate-optimal'' behavior with a broad effective range (roughly 30\%--80\%), suggesting that the framework is not overly sensitive to precise DR tuning.
\end{document}